\documentclass{amsproc}
\usepackage{amssymb}
\usepackage{amsmath}
\usepackage{eurosym}
\usepackage{amsfonts}

            \advance \topmargin by -\headheight \advance
            \topmargin by -\headsep
            \evensidemargin
            \oddsidemargin
\newtheorem{theorem}{Theorem}[section]
\newtheorem{lemma}[theorem]{Lemma}
\newtheorem{definition}[theorem]{Definition}
\newtheorem{corollary}[theorem]{Corollary}

\numberwithin{equation}{section}
\input{tcilatex}

\begin{document}
\title[Tribonacci and Tribonacci-Lucas Matrix Sequences with Negative
Subscripts]{Tribonacci and Tribonacci-Lucas Matrix Sequences with Negative
Subscripts}
\thanks{}
\author[Y\"{u}ksel ~Soykan]{Y\"{U}KSEL SOYKAN}
\maketitle

\begin{center}
\textsl{Zonguldak B\"{u}lent Ecevit University, Department of Mathematics, }

\textsl{Art and Science Faculty, 67100, Zonguldak, Turkey }

\textsl{e-mail: \ yuksel\_soykan@hotmail.com}
\end{center}

\textbf{Abstract.} In this paper, we define Tribonacci and Tribonacci-Lucas
matrix sequences with negative indices and investigate their properties.

\textbf{2010 Mathematics Subject Classification.} 11B39, 11B83.

\textbf{Keywords. }Tribonacci numbers, Tribonacci matrix sequence,
Tribonacci-Lucas matrix sequence.

\section{Introduction and Preliminaries}

Fibonacci sequence is often used as a model of recursive phenomena in
physics and engineering (see for example, [\ref{akbulak2009}]), chemistry
(see for instance, [\ref{randic1996}], [\ref{yegnanarayanan2012}]), botany
(see e.g. [\ref{ridley1982}]), medicine (see e.g. [\ref{cohen2002}]). As
Fibonacci sequence, also Tribonacci sequence has many applications to such
as coding theory (see e.g. [\ref{basu2014}]), game theory (see e.g. [\ref%
{duchene2008}] and references therein).

By a Mathematical point of view, there are analogies properties between
Horadam sequences (such as Fibonacci and Lucas sequences) and generalized
Tribonacci sequences (such as Padovan and usual Tribonacci sequences). For
those concerning asymptotic process if we look at the characteristic
polynomial $x^{2}-x-1=0$ associated to Fibonacci recursive relation, we have
that it has a unique (real) root of maximum modulus, that is also the limit
of the ratio of two consecutive Fibonacci numbers:%
\begin{equation*}
\lim_{n\rightarrow \infty }\frac{F_{n+1}}{F_{n}}=\alpha _{F}
\end{equation*}%
where $\alpha _{F}:=\frac{1+\sqrt{5}}{2}$ denotes the highly celebrated
Golden mean (also called Golden section or Golden ratio), (see Example 4.1
in [\ref{fiorenza2011}]).

Similarly, if we consider the characteristic polynomial $x^{3}-x-1=0$
associated to Padovan (a generalized Tribonacci ) recursive relation, we
have that it has a unique (real) root of maximum modulus, that is also the
limit of the ratio of two consecutive Padovan numbers:%
\begin{equation*}
\lim_{n\rightarrow \infty }\frac{P_{n+1}}{P_{n}}=\alpha _{P}
\end{equation*}%
where $\alpha _{P}:=\sqrt[3]{\frac{1}{2}+\frac{1}{6}\sqrt{\frac{23}{3}}}+%
\sqrt[3]{\frac{1}{2}-\frac{1}{6}\sqrt{\frac{23}{3}}}$ denotes the Plastic
ratio (see Example 4.9 in [\ref{fiorenza2011}]) which have many applications
to such as architecture, see [\ref{bib:marohnic2012}]).

In our case, if we consider the characteristic polynomial $x^{3}-x^{2}-x-1=0$
associated to Tribonacci recursive relation, we have that it has a unique
(real) root of maximum modulus, that is also the limit of the ratio of two
consecutive Tribonacci numbers:%
\begin{equation*}
\lim_{n\rightarrow \infty }\frac{T_{n+1}}{T_{n}}=\alpha _{T}
\end{equation*}%
where $\alpha _{T}:=\frac{1+\sqrt[3]{19+3\sqrt{33}}+\sqrt[3]{19-3\sqrt{33}}}{%
3}$ denotes the Tribonacci ratio (see for example Example 4.9 in [\ref%
{fiorenza2011}] or for a basic proof see [\ref{bueno:2015}]). For a short
introduction to these three constants, see [\ref{piezas}].

In fact, for linear homogeneous recursive sequences with constant
coefficients, Fiorenza and Vincenzi [\ref{fiorenza2011}] find a necessary
and sufficient condition for the existence of the limit of the ratio of
consecutive terms. As a corollary of their results, the limit of the ratio
of adjacent terms is characterized as the unique leading root of the
characteristic polynomial.

On the other hand, the matrix sequences have taken so much interest for
different type of numbers. For matrix sequences of generalized Horadam type
numbers, see for example [\ref{civciv2008}], [\ref{civciv2008b}], [\ref%
{gulec2012}], [\ref{uslu2013}], [\ref{uygun2015}], [\ref{uygun2016}], [\ref%
{yazlik2013a}], [\ref{wani2018}], and for matrix sequences of generalized
Tribonacci type numbers, see for instance [\ref{cerdamoralez2018aa}], [\ref%
{yilmaz2013}], [\ref{yilmaz2014a}],\ [\ref{soykan2018triposmatrix}].

In this paper, for negative indices, the matrix sequences of Tribonacci and
Tribonacci-Lucas numbers will be defined. Then, by giving the generating
functions, the Binet formulas, and summation formulas over these new matrix
sequences, we will obtain some fundamental properties on Tribonacci and
Tribonacci-Lucas numbers. Also, we will present the relationship between
these matrix sequences.

Now, we give some background about Tribonacci and Tribonacci-Lucas numbers.
Tribonacci sequence $\{T_{n}\}_{n\geq 0}$ (sequence A000073 in [\ref%
{bib:sloane}]) and Tribonacci-Lucas sequence $\{K_{n}\}_{n\geq 0}$ (sequence
A001644 in [\ref{bib:sloane}]) are defined by the third-order recurrence
relations 
\begin{equation}
T_{n}=T_{n-1}+T_{n-2}+T_{n-3},\text{ \ \ \ \ }T_{0}=0,T_{1}=1,T_{2}=1,
\label{equati:fvcvxghsbnz}
\end{equation}%
and 
\begin{equation}
K_{n}=K_{n-1}+K_{n-2}+K_{n-3},\text{ \ \ \ \ }K_{0}=3,K_{1}=1,K_{2}=3
\label{equati:pazertvbcunsmn}
\end{equation}%
respectively. Tribonacci concept was introduced by M. Feinberg\ [\ref%
{bib:feinberg1963}] in 1963. Basic properties of it is given in [\ref%
{bib:bruce1984}], [\ref{bib:scott1977}], [\ref{bib:shannon1977}], [\ref%
{bib:spickerman1981}], [\ref{bib:yalavigi1972}].

The sequences $\{T_{n}\}_{n\geq 0}$ and $\{K_{n}\}_{n\geq 0}$ can be
extended to negative subscripts by defining%
\begin{equation*}
T_{-n}=-T_{-(n-1)}-T_{-(n-2)}+T_{-(n-3)}
\end{equation*}%
and%
\begin{equation*}
K_{-n}=-K_{-(n-1)}-K_{-(n-2)}+K_{-(n-3)}
\end{equation*}%
for $n=1,2,3,...$ respectively. Therefore, recurrences (\ref%
{equati:fvcvxghsbnz}) and (\ref{equati:pazertvbcunsmn}) hold for all integer 
$n.$

We can give some relations between $\{T_{n}\}$ and $\{K_{n}\}$\ as 
\begin{equation}
K_{n}=3T_{n+1}-2T_{n}-T_{n-1}  \label{equation:gfdvbxczvsadou}
\end{equation}%
and%
\begin{equation}
K_{n}=T_{n}+2T_{n-1}+3T_{n-2}  \label{equation:fvcdxsartewqa}
\end{equation}%
and also%
\begin{equation}
K_{n}=4T_{n+1}-T_{n}-T_{n+2}.  \label{equation:yuhtgsdafcscvb}
\end{equation}%
Note that the last three identities hold for all integers $n.$

The first few Tribonacci numbers and Tribonacci Lucas numbers with positive
subscript are given in the following table:

$%
\begin{array}{ccccccccccccccc}
n & 0 & 1 & 2 & 3 & 4 & 5 & 6 & 7 & 8 & 9 & 10 & 11 & 12 & ... \\ 
T_{n} & 0 & 1 & 1 & 2 & 4 & 7 & 13 & 24 & 44 & 81 & 149 & 274 & 504 & ... \\ 
T_{-n} & 0 & 0 & 1 & -1 & 0 & 2 & -3 & 1 & 4 & -8 & 5 & 7 & -20 & ...%
\end{array}%
$

The first few Tribonacci numbers and Tribonacci Lucas numbers with negative
subscript are given in the following table:

$%
\begin{array}{ccccccccccccccc}
n & 0 & 1 & 2 & 3 & 4 & 5 & 6 & 7 & 8 & 9 & 10 & 11 & 12 & ... \\ 
K_{n} & 3 & 1 & 3 & 7 & 11 & 21 & 39 & 71 & 131 & 241 & 443 & 815 & 1499 & 
... \\ 
K_{-n} & 3 & -1 & -1 & 5 & -5 & -1 & 11 & -15 & 3 & 23 & -41 & 21 & 43 & ...%
\end{array}%
$

It is well known that for all integers $n,$ usual Tribonacci and
Tribonacci-Lucas numbers can be expressed using Binet's formulas%
\begin{equation}
T_{n}=\frac{\alpha ^{n+1}}{(\alpha -\beta )(\alpha -\gamma )}+\frac{\beta
^{n+1}}{(\beta -\alpha )(\beta -\gamma )}+\frac{\gamma ^{n+1}}{(\gamma
-\alpha )(\gamma -\beta )}  \label{equat:mnopcvbedcxzsa}
\end{equation}%
and%
\begin{equation}
K_{n}=\alpha ^{n}+\beta ^{n}+\gamma ^{n}  \label{equation:cfrdcsxszouea}
\end{equation}%
respectively, where $\alpha ,\beta $ and $\gamma $ are the roots of the
cubic equation $x^{3}-x^{2}-x-1=0.$ Moreover, 
\begin{eqnarray*}
\alpha &=&\frac{1+\sqrt[3]{19+3\sqrt{33}}+\sqrt[3]{19-3\sqrt{33}}}{3}, \\
\beta &=&\frac{1+\omega \sqrt[3]{19+3\sqrt{33}}+\omega ^{2}\sqrt[3]{19-3%
\sqrt{33}}}{3}, \\
\gamma &=&\frac{1+\omega ^{2}\sqrt[3]{19+3\sqrt{33}}+\omega \sqrt[3]{19-3%
\sqrt{33}}}{3}
\end{eqnarray*}%
where%
\begin{equation*}
\omega =\frac{-1+i\sqrt{3}}{2}=\exp (2\pi i/3),
\end{equation*}%
is a primitive cube root of unity. Note that we have the following identities%
\begin{eqnarray*}
\alpha +\beta +\gamma &=&1, \\
\alpha \beta +\alpha \gamma +\beta \gamma &=&-1, \\
\alpha \beta \gamma &=&1.
\end{eqnarray*}

The generating functions for the Tribonacci sequence $\{T_{n}\}_{n\geq 0}$
and Tribonacci-Lucas sequence $\{K_{n}\}_{n\geq 0}$ are%
\begin{equation}
\sum_{n=0}^{\infty }T_{n}x^{n}=\frac{x}{1-x-x^{2}-x^{3}}\text{ \ and \ }%
\sum_{n=0}^{\infty }K_{n}x^{n}=\frac{3-2x-x^{2}}{1-x-x^{2}-x^{3}}.
\label{equation:yugdfvxbgsopqac}
\end{equation}

Note that the Binet form of a sequence satisfying (\ref{equati:fvcvxghsbnz}%
)\ and (\ref{equati:pazertvbcunsmn}) for non-negative integers is valid for
all integers $n.$ This result of Howard and Saidak [\ref{bib:howard2010}] is
even true in the case of higher-order recurrence relations as the following
theorem shows.

\begin{theorem}[{[\protect\ref{bib:howard2010}]}]
\label{theorem:fvgxdfsxczsaer}Let $\{w_{n}\}$ be a sequence such that 
\begin{equation*}
\{w_{n}\}=a_{1}w_{n-1}+a_{2}w_{n-2}+...+a_{k}w_{n-k}
\end{equation*}%
for all integers $n,$ with arbitrary initial conditions $%
w_{0},w_{1},...,w_{k-1}.$ Assume that each $a_{i}$ and the initial
conditions are complex numbers. Write%
\begin{eqnarray}
f(x) &=&x^{k}-a_{1}x^{k-1}-a_{2}x^{k-2}-...-a_{k-1}x-a_{k}
\label{equation:mnbvyuhgoewapvbc} \\
&=&(x-\alpha _{1})^{d_{1}}(x-\alpha _{2})^{d_{2}}...(x-\alpha _{h})^{d_{h}} 
\notag
\end{eqnarray}%
with $d_{1}+d_{2}+...+d_{h}=k,$ and $\alpha _{1},\alpha _{2},...,\alpha _{k}$
distinct. Then

\begin{description}
\item[(a)] For all $n,$%
\begin{equation}
w_{n}=\sum_{m=1}^{k}N(n,m)(\alpha _{m})^{n}
\label{equation:yuhnbvdfscxzratqw}
\end{equation}%
where%
\begin{equation*}
N(n,m)=A_{1}^{(m)}+A_{2}^{(m)}n+...+A_{r_{m}}^{(m)}n^{r_{m}-1}=%
\sum_{u=0}^{r_{m}-1}A_{u+1}^{(m)}n^{u}
\end{equation*}%
with each $A_{i}^{(m)}$ a constant determined by the initial conditions for $%
\{w_{n}\}$. Here, equation (\ref{equation:yuhnbvdfscxzratqw}) is called the
Binet form (or Binet formula) for $\{w_{n}\}.$ We assume that $f(0)\neq 0$
so that $\{w_{n}\}$ can be extended to negative integers $n.$

If the zeros of (\ref{equation:mnbvyuhgoewapvbc}) are distinct, as they are
in our examples, then%
\begin{equation*}
w_{n}=A_{1}(\alpha _{1})^{n}+A_{2}(\alpha _{2})^{n}+...+A_{k}(\alpha
_{k})^{n}.
\end{equation*}

\item[(b)] The Binet form for $\{w_{n}\}$ is valid for all integers $n.$
\end{description}
\end{theorem}

In [\ref{soykan2018triposmatrix}], Soykan introduced the following
definition of Tribonacci and Tribonacci-Lucas matrix sequences and
investigated their properties.

\begin{definition}
For any integer $n\geq 0,$ the Tribonacci matrix $(\mathcal{T}_{n})$ and
Tribonacci-Lucas matrix $(\mathcal{K}_{n})$ are defined by%
\begin{eqnarray}
\mathcal{T}_{n} &=&\mathcal{T}_{n-1}+\mathcal{T}_{n-2}+\mathcal{T}_{n-3},
\label{equation:yusoabpsmnbscv} \\
\mathcal{K}_{n} &=&\mathcal{K}_{n-1}+\mathcal{K}_{n-2}+\mathcal{K}_{n-3,}
\label{equation:dfscxzvayuewrtsfg}
\end{eqnarray}%
respectively, with initial conditions%
\begin{equation*}
\mathcal{T}_{0}=\left( 
\begin{array}{ccc}
1 & 0 & 0 \\ 
0 & 1 & 0 \\ 
0 & 0 & 1%
\end{array}%
\right) ,\mathcal{T}_{1}=\left( 
\begin{array}{ccc}
1 & 1 & 1 \\ 
1 & 0 & 0 \\ 
0 & 1 & 0%
\end{array}%
\right) ,\mathcal{T}_{2}=\left( 
\begin{array}{ccc}
2 & 2 & 1 \\ 
1 & 1 & 1 \\ 
1 & 0 & 0%
\end{array}%
\right)
\end{equation*}%
and%
\begin{equation*}
\mathcal{K}_{0}=\left( 
\begin{array}{ccc}
1 & 2 & 3 \\ 
3 & -2 & -1 \\ 
-1 & 4 & -1%
\end{array}%
\right) ,\mathcal{K}_{1}=\left( 
\begin{array}{ccc}
3 & 4 & 1 \\ 
1 & 2 & 3 \\ 
3 & -2 & -1%
\end{array}%
\right) ,\mathcal{K}_{2}=\left( 
\begin{array}{ccc}
7 & 4 & 3 \\ 
3 & 4 & 1 \\ 
1 & 2 & 3%
\end{array}%
\right) .
\end{equation*}
\end{definition}

$(\mathcal{T}_{n})$ and $(\mathcal{K}_{n})$ have some good properties which
is given in next two Theorems.

\begin{theorem}[{[\protect\ref{soykan2018triposmatrix}]}]
\label{theorem:olscdxfczdarepwout}For all non-negative integers $m$ and $n,\ 
$we have the following identities.

\begin{description}
\item[(a)] $\mathcal{T}_{m}\mathcal{T}_{n}=\mathcal{T}_{m+n}=\mathcal{T}_{n}%
\mathcal{T}_{m},$

\item[(b)] $\mathcal{T}_{m}\mathcal{K}_{n}=\mathcal{K}_{n}\mathcal{T}_{m}=%
\mathcal{K}_{m+n},$

\item[(c)] $\mathcal{K}_{m}\mathcal{K}_{n}=\mathcal{K}_{n}\mathcal{K}_{m}=9%
\mathcal{T}_{m+n+2}-12\mathcal{T}_{m+n+1}+\mathcal{T}_{m+n}+\mathcal{T}%
_{m+n-1}+\mathcal{T}_{m+n-2},$

\item[(d)] $\mathcal{K}_{m}\mathcal{K}_{n}=\mathcal{K}_{n}\mathcal{K}_{m}=%
\mathcal{T}_{m+n}+4\mathcal{T}_{m+n-1}+10\mathcal{T}_{m+n-2}+12\mathcal{T}%
_{m+n-3}+\allowbreak 9\mathcal{T}_{m+n-4},$

\item[(e)] $\mathcal{K}_{m}\mathcal{K}_{n}=\mathcal{K}_{n}\mathcal{K}_{m}=%
\mathcal{T}_{m+n}-8\mathcal{T}_{m+n+1}+18\mathcal{T}_{m+n+2}-8\mathcal{T}%
_{m+n+3}+\mathcal{T}_{m+n+4}.$
\end{description}
\end{theorem}

We now give the Binet formulas for the Tribonacci and Tribonacci-Lucas
matrix sequences.

\begin{theorem}[{[\protect\ref{soykan2018triposmatrix}]}]
\label{theorem:sacasdmnbhgvc}For every integer $n,$ the Binet formulas of
the Tribonacci and Tribonacci-Lucas matrix sequences are given by%
\begin{eqnarray}
\mathcal{T}_{n} &=&A_{1}\alpha ^{n}+B_{1}\beta ^{n}+C_{1}\gamma ^{n},
\label{equation:mousxzadfgytc} \\
\mathcal{K}_{n} &=&A_{2}\alpha ^{n}+B_{2}\beta ^{n}+C_{2}\gamma ^{n}.
\label{equat:fvbxdszuqwsazx}
\end{eqnarray}%
where%
\begin{eqnarray*}
A_{1} &=&\frac{\alpha \mathcal{T}_{2}+\alpha (\alpha -1)\mathcal{T}_{1}+%
\mathcal{T}_{0}}{\alpha \left( \alpha -\gamma \right) \left( \alpha -\beta
\right) },B_{1}=\frac{\beta \mathcal{T}_{2}+\beta (\beta -1)\mathcal{T}_{1}+%
\mathcal{T}_{0}}{\beta \left( \beta -\gamma \right) \left( \beta -\alpha
\right) },C_{1}=\frac{\gamma \mathcal{T}_{2}+\gamma (\gamma -1)\mathcal{T}%
_{1}+\mathcal{T}_{0}}{\gamma \left( \gamma -\beta \right) \left( \gamma
-\alpha \right) } \\
A_{2} &=&\frac{\alpha \mathcal{K}_{2}+\alpha (\alpha -1)\mathcal{K}_{1}+%
\mathcal{K}_{0}}{\alpha \left( \alpha -\gamma \right) \left( \alpha -\beta
\right) },B_{2}=\frac{\beta \mathcal{K}_{2}+\beta (\beta -1)\mathcal{K}_{1}+%
\mathcal{K}_{0}}{\beta \left( \beta -\gamma \right) \left( \beta -\alpha
\right) },C_{2}=\frac{\gamma \mathcal{K}_{2}+\gamma (\gamma -1)\mathcal{K}%
_{1}+\mathcal{K}_{0}}{\gamma \left( \gamma -\beta \right) \left( \gamma
-\alpha \right) }.
\end{eqnarray*}
\end{theorem}

Note that the Binet formulas given above hold for all integers $n.$

\section{The Matrix Sequences of Negative Subscripts Tribonacci and
Tribonacci-Lucas Numbers}

The sequences $\{\mathcal{T}_{n}\}_{n\geq 0}$ and $\{\mathcal{K}%
_{n}\}_{n\geq 0}$ can be extended to negative subscripts by defining 
\begin{equation*}
\mathcal{T}_{-n}=-\mathcal{T}_{-(n-1)}-\mathcal{T}_{-(n-2)}+\mathcal{T}%
_{-(n-3)}
\end{equation*}%
and%
\begin{equation*}
\mathcal{K}_{-n}=-\mathcal{K}_{-(n-1)}-\mathcal{K}_{-(n-2)}+\mathcal{K}%
_{-(n-3)}
\end{equation*}%
for $n=1,2,3,...$ respectively. Therefore, recurrences (\ref%
{equation:yusoabpsmnbscv}) and (\ref{equation:dfscxzvayuewrtsfg}) hold for
all integer. i.e. starting with $n=-1$ and working backwards, we extend $\{%
\mathcal{T}_{n}\}$ and $\{\mathcal{K}_{n}\}$ to negative indices. The first
few Tribonacci numbers and Tribonacci Lucas numbers with negative subscript
can be found as follows:%
\begin{eqnarray*}
\mathcal{T}_{-1} &=&\mathcal{T}_{2}-\mathcal{T}_{1}-\mathcal{T}_{0}=\left( 
\begin{array}{ccc}
0 & 1 & 0 \\ 
0 & 0 & 1 \\ 
1 & -1 & -1%
\end{array}%
\right) \\
\mathcal{T}_{-2} &=&\mathcal{T}_{1}-\mathcal{T}_{0}-\mathcal{T}_{-1}=\left( 
\begin{array}{ccc}
0 & 0 & 1 \\ 
1 & -1 & -1 \\ 
-1 & 2 & 0%
\end{array}%
\right) \\
\mathcal{T}_{-3} &=&\mathcal{T}_{0}-\mathcal{T}_{-1}-\mathcal{T}_{-2}=\left( 
\begin{array}{ccc}
1 & -1 & -1 \\ 
-1 & 2 & 0 \\ 
0 & -1 & 2%
\end{array}%
\right)
\end{eqnarray*}%
and%
\begin{eqnarray*}
\mathcal{K}_{-1} &=&\mathcal{K}_{2}-\mathcal{K}_{1}-\mathcal{K}_{0}=\left( 
\begin{array}{ccc}
3 & -2 & -1 \\ 
-1 & 4 & -1 \\ 
-1 & 0 & 5%
\end{array}%
\right) \\
\mathcal{K}_{-2} &=&\mathcal{K}_{1}-\mathcal{K}_{0}-\mathcal{K}_{-1}=\left( 
\begin{array}{ccc}
-1 & 4 & -1 \\ 
-1 & 0 & 5 \\ 
5 & -6 & -5%
\end{array}%
\right) \\
\mathcal{K}_{-3} &=&\mathcal{K}_{0}-\mathcal{K}_{-1}-\mathcal{K}_{-2}=\left( 
\begin{array}{ccc}
-1 & 0 & 5 \\ 
5 & -6 & -5 \\ 
-5 & 10 & -1%
\end{array}%
\right) .
\end{eqnarray*}%
Actually, we can formally define $\{\mathcal{T}_{n}\}$ and $\{\mathcal{K}%
_{n}\}$ for negative indices as follows.

\begin{definition}
\label{definition::negabvhsyteur}For any integer $n\geq 0,$ the negative
indices Tribonacci matrix $(\mathcal{T}_{-n})$ and Tribonacci-Lucas matrix $(%
\mathcal{K}_{-n})$ are defined by%
\begin{equation}
\mathcal{T}_{-n}=\mathcal{T}_{-(n-3)}-\mathcal{T}_{-(n-1)}-\mathcal{T}%
_{-(n-2)}=\mathcal{T}_{-n+3}-\mathcal{T}_{-n+2}-\mathcal{T}_{-n+1}
\label{equation:posmbnkstydhg}
\end{equation}%
and%
\begin{equation}
\mathcal{K}_{-n}=\mathcal{K}_{-(n-3)}-\mathcal{K}_{-(n-1)}-\mathcal{K}%
_{-(n-2)}=\mathcal{K}_{-n+3}-\mathcal{K}_{-n+2}-\mathcal{K}_{-n+1}
\end{equation}%
respectively, with initial conditions%
\begin{equation}
\mathcal{T}_{0}=\left( 
\begin{array}{ccc}
1 & 0 & 0 \\ 
0 & 1 & 0 \\ 
0 & 0 & 1%
\end{array}%
\right) ,\mathcal{T}_{-1}=\left( 
\begin{array}{ccc}
0 & 1 & 0 \\ 
0 & 0 & 1 \\ 
1 & -1 & -1%
\end{array}%
\right) ,\mathcal{T}_{-2}=\left( 
\begin{array}{ccc}
0 & 0 & 1 \\ 
1 & -1 & -1 \\ 
-1 & 2 & 0%
\end{array}%
\right)  \label{equati:erkasnbchdutrd}
\end{equation}%
and%
\begin{equation}
\mathcal{K}_{0}=\left( 
\begin{array}{ccc}
1 & 2 & 3 \\ 
3 & -2 & -1 \\ 
-1 & 4 & -1%
\end{array}%
\right) ,\mathcal{K}_{-1}=\left( 
\begin{array}{ccc}
3 & -2 & -1 \\ 
-1 & 4 & -1 \\ 
-1 & 0 & 5%
\end{array}%
\right) ,\mathcal{K}_{-2}=\left( 
\begin{array}{ccc}
-1 & 4 & -1 \\ 
-1 & 0 & 5 \\ 
5 & -6 & -5%
\end{array}%
\right) .  \label{equat:negbvcagosetrya}
\end{equation}
\end{definition}

\ The following theorem gives the $n$th general terms of the Tribonacci and
Tribonacci-Lucas matrix sequences with negative indices.

\begin{theorem}
\label{theorem:gbvdcsxfazoerstaqw}For any integer $n\geq 0,$ we have the
following formulas of the matrix sequences:%
\begin{eqnarray}
\mathcal{T}_{-n} &=&\left( 
\begin{array}{ccc}
T_{-n+1} & T_{-n}+T_{-n-1} & T_{-n} \\ 
T_{-n} & T_{-n-1}+T_{-n-2} & T_{-n-1} \\ 
T_{-n-1} & T_{-n-2}+T_{-n-3} & T_{-n-2}%
\end{array}%
\right)  \label{equation:dfgvcoerdwfacxcv} \\
\mathcal{K}_{-n} &=&\left( 
\begin{array}{ccc}
K_{-n+1} & K_{-n}+K_{-n-1} & K_{-n} \\ 
K_{-n} & K_{-n-1}+K_{-n-2} & K_{-n-1} \\ 
K_{-n-1} & K_{-n-2}+K_{-n-3} & K_{-n-2}%
\end{array}%
\right) .  \label{equation:dcvbosavzbxyuredfwqs}
\end{eqnarray}
\end{theorem}

Proof. We prove (\ref{equation:dfgvcoerdwfacxcv}) by strong mathematical
induction on $n$. (\ref{equation:dcvbosavzbxyuredfwqs}) can be proved
similarly.

If $n=0$ then since $T_{1}=1,T_{0}=T_{-1}=0$,\ $T_{-2}=1,$ $%
T_{-3}=-1,T_{-4}=0,$ we have 
\begin{equation*}
\mathcal{T}_{0}=\left( 
\begin{array}{ccc}
T_{1} & T_{0}+T_{-1} & T_{0} \\ 
T_{0} & T_{-1}+T_{-2} & T_{-1} \\ 
T_{-1} & T_{-2}+T_{-3} & T_{-2}%
\end{array}%
\right) =\left( 
\begin{array}{ccc}
1 & 0 & 0 \\ 
0 & 1 & 0 \\ 
0 & 0 & 1%
\end{array}%
\right)
\end{equation*}%
which is true and 
\begin{equation*}
\mathcal{T}_{-1}=\left( 
\begin{array}{ccc}
T_{0} & T_{-1}+T_{-2} & T_{-1} \\ 
T_{-1} & T_{-2}+T_{-3} & T_{-2} \\ 
T_{-2} & T_{-3}+T_{-4} & T_{-3}%
\end{array}%
\right) =\left( 
\begin{array}{ccc}
0 & 1 & 0 \\ 
0 & 0 & 1 \\ 
1 & -1 & -1%
\end{array}%
\right)
\end{equation*}%
which is true. Assume that the equality holds for $n\leq k.$ For $n=k+1,$ we
have%
\begin{eqnarray*}
\mathcal{T}_{-(k+1)} &=&\mathcal{T}_{-(k+1)+3}-\mathcal{T}_{-(k+1)+2}-%
\mathcal{T}_{-(k+1)+1}=\mathcal{T}_{-(k-2)}-\mathcal{T}_{-(k-1)}-\mathcal{T}%
_{-k} \\
&=&\mathcal{T}_{2-k}-\mathcal{T}_{1-k}-\mathcal{T}_{-k} \\
&=&\left( 
\begin{array}{ccc}
T_{3-k} & T_{1-k}+T_{2-k} & T_{2-k} \\ 
T_{2-k} & T_{-k}+T_{1-k} & T_{1-k} \\ 
T_{1-k} & T_{-k}+T_{-k-1} & T_{-k}%
\end{array}%
\right) -\left( 
\begin{array}{ccc}
T_{2-k} & T_{-k}+T_{1-k} & T_{1-k} \\ 
T_{1-k} & T_{-k}+T_{-k-1} & T_{-k} \\ 
T_{-k} & T_{-k-1}+T_{-k-2} & T_{-k-1}%
\end{array}%
\right) \\
&&-\left( 
\begin{array}{ccc}
T_{1-k} & T_{-k}+T_{-k-1} & T_{-k} \\ 
T_{-k} & T_{-k-1}+T_{-k-2} & T_{-k-1} \\ 
T_{-k-1} & T_{-k-2}+T_{-k-3} & T_{-k-2}%
\end{array}%
\right) \\
&=&\left( 
\begin{array}{ccc}
T_{-k} & T_{-k-1}+T_{-k-2} & T_{-k-1} \\ 
T_{-k-1} & T_{-k-2}+T_{-k-3} & T_{-k-2} \\ 
T_{-k-2} & T_{-k-3}+T_{-k-4} & T_{-k-3}%
\end{array}%
\right) \\
&=&\left( 
\begin{array}{ccc}
T_{-(k+1)+1} & T_{-(k+1)}+T_{-(k+1)-1} & T_{-(k+1)} \\ 
T_{-(k+1)} & T_{-(k+1)-1}+T_{-(k+1)-2} & T_{-(k+1)-1} \\ 
T_{-(k+1)-1} & T_{-(k+1)-2}+T_{-(k+1)-3} & T_{-(k+1)-2}%
\end{array}%
\right)
\end{eqnarray*}%
Thus, by strong induction on $n,$ this proves (\ref%
{equation:dfgvcoerdwfacxcv}). 
\endproof%

The Binet formulas for the matrix sequences of Tribonacci and
Tribonacci-Lucas numbers are given in [\ref{soykan2018triposmatrix}], see
Theorem \ref{theorem:sacasdmnbhgvc} above. For the completeness of the
paper, we now give the Binet formula for the Tribonacci and Tribonacci-Lucas
matrix sequences with negative indices.

\begin{theorem}
\label{theorem:sdfvbaxdreqwsax}For every non-negative integer $n,$ the Binet
formulas of the Tribonacci and Tribonacci-Lucas matrix sequences are given by%
\begin{eqnarray}
\mathcal{T}_{-n} &=&A_{3}\alpha ^{-n}+B_{3}\beta ^{-n}+C_{3}\gamma ^{-n},
\label{equati:ghnvcxdtsrtyuaxz} \\
\mathcal{K}_{-n} &=&A_{4}\alpha ^{-n}+B_{4}\beta ^{-n}+C_{4}\gamma ^{-n}.
\label{equat:nbfdbvcxsdfzusrat}
\end{eqnarray}%
where%
\begin{eqnarray*}
A_{3} &=&\frac{\alpha \mathcal{T}_{-2}+(\alpha -1)\alpha ^{2}\mathcal{T}%
_{-1}+\alpha ^{2}\mathcal{T}_{0}}{\left( \alpha -\gamma \right) \left(
\alpha -\beta \right) },B_{3}=\frac{\beta \mathcal{T}_{-2}+(\beta -1)\beta
^{2}\mathcal{T}_{-1}+\beta ^{2}\mathcal{T}_{0}}{\left( \beta -\gamma \right)
\left( \beta -\alpha \right) },C_{3}=\frac{\gamma \mathcal{T}_{-2}+(\gamma
-1)\gamma ^{2}\mathcal{T}_{-1}+\gamma ^{2}\mathcal{T}_{0}}{\left( \gamma
-\beta \right) \left( \gamma -\alpha \right) } \\
A_{4} &=&\frac{\alpha \mathcal{K}_{-2}+(\alpha -1)\alpha ^{2}\mathcal{K}%
_{-1}+\alpha ^{2}\mathcal{K}_{0}}{\left( \alpha -\gamma \right) \left(
\alpha -\beta \right) },B_{4}=\frac{\beta \mathcal{K}_{-2}+(\beta -1)\beta
^{2}\mathcal{K}_{-1}+\beta ^{2}\mathcal{K}_{0}}{\left( \beta -\gamma \right)
\left( \beta -\alpha \right) },C_{4}=\frac{\gamma \mathcal{K}_{-2}+(\gamma
-1)\gamma ^{2}\mathcal{K}_{-1}+\gamma ^{2}\mathcal{K}_{0}}{\left( \gamma
-\beta \right) \left( \gamma -\alpha \right) }.
\end{eqnarray*}
\end{theorem}

\textit{Proof.} Note that the proof is based on the recurrence relations (%
\ref{equati:erkasnbchdutrd}) and (\ref{equat:negbvcagosetrya}) in Definition %
\ref{definition::negabvhsyteur}.

We prove (\ref{equati:ghnvcxdtsrtyuaxz}). By the assumption, the
characteristic equation of (\ref{equation:posmbnkstydhg}) is $%
x^{3}+x^{2}+x-1=0$ and the roots of it are $\frac{1}{\alpha },\frac{1}{\beta 
}$ and $\frac{1}{\gamma }.$ So it's general solution is given by%
\begin{equation*}
\mathcal{T}_{-n}=A_{3}\alpha ^{-n}+B_{3}\beta ^{-n}+C_{3}\gamma ^{-n}.
\end{equation*}%
Using initial condition which is given in Definition \ref%
{definition::negabvhsyteur}, and also applying lineer algebra operations, we
obtain the matrices $A_{3},B_{3},C_{3}$ as desired. This gives the formula
for $\mathcal{T}_{-n}.$

Similarly we have the formula (\ref{equat:nbfdbvcxsdfzusrat}). 
\endproof%

In fact, again by Theorem \ref{theorem:fvgxdfsxczsaer}, Theorem \ref%
{theorem:sdfvbaxdreqwsax} is true for all integers $n.$ If we compare
Theorem \ref{theorem:sacasdmnbhgvc} and Theorem \ref{theorem:sdfvbaxdreqwsax}%
\ we obtain%
\begin{equation}
A_{1}=A_{3},\text{ }B_{1}=B_{3},\text{ }C_{1}=C_{3},
\label{equation:hbvcdfrgsxzasewq}
\end{equation}%
and 
\begin{equation}
A_{2}=A_{4},\text{ }B_{2}=B_{4},\text{ }C_{2}=C_{4},
\end{equation}%
i.e., 
\begin{eqnarray*}
\alpha \mathcal{T}_{-2}+(\alpha -1)\alpha ^{2}\mathcal{T}_{-1}+\alpha ^{2}%
\mathcal{T}_{0} &=&\frac{\alpha \mathcal{T}_{2}+\alpha (\alpha -1)\mathcal{T}%
_{1}+\mathcal{T}_{0}}{\alpha }, \\
\beta \mathcal{T}_{-2}+(\beta -1)\beta ^{2}\mathcal{T}_{-1}+\beta ^{2}%
\mathcal{T}_{0} &=&\frac{\beta \mathcal{T}_{2}+\beta (\beta -1)\mathcal{T}%
_{1}+\mathcal{T}_{0}}{\beta }, \\
\gamma \mathcal{T}_{-2}+(\gamma -1)\gamma ^{2}\mathcal{T}_{-1}+\gamma ^{2}%
\mathcal{T}_{0} &=&\frac{\gamma \mathcal{T}_{2}+\gamma (\gamma -1)\mathcal{T}%
_{1}+\mathcal{T}_{0}}{\gamma }, \\
\alpha \mathcal{K}_{-2}+(\alpha -1)\alpha ^{2}\mathcal{K}_{-1}+\alpha ^{2}%
\mathcal{K}_{0} &=&\frac{\alpha \mathcal{K}_{2}+\alpha (\alpha -1)\mathcal{K}%
_{1}+\mathcal{K}_{0}}{\alpha }, \\
\beta \mathcal{K}_{-2}+(\beta -1)\beta ^{2}\mathcal{K}_{-1}+\beta ^{2}%
\mathcal{K}_{0} &=&\frac{\beta \mathcal{K}_{2}+\beta (\beta -1)\mathcal{K}%
_{1}+\mathcal{K}_{0}}{\beta }, \\
\gamma \mathcal{K}_{-2}+(\gamma -1)\gamma ^{2}\mathcal{K}_{-1}+\gamma ^{2}%
\mathcal{K}_{0} &=&\frac{\gamma \mathcal{K}_{2}+\gamma (\gamma -1)\mathcal{K}%
_{1}+\mathcal{K}_{0}}{\gamma }.
\end{eqnarray*}

The well known Binet formulas (for positive and negative indices) for
Tribonacci and Tribonacci-Lucas numbers are given in (\ref%
{equat:mnopcvbedcxzsa})\ and (\ref{equation:cfrdcsxszouea}) respectively.
But, for negative indices, we will obtain these functions in terms of
Tribonacci and Tribonacci-Lucas matrix sequences as a consequence of
Theorems \ref{theorem:gbvdcsxfazoerstaqw} and \ref{theorem:sdfvbaxdreqwsax}.
To do this, we will give the formulas for these numbers by means of the
related matrix sequences. In fact, in the proof of next corollary, we will
just compare the linear combination of the 2nd row and 1st column entries of
the matrices.

\begin{corollary}
\label{corollary:fvxgfsypqavbxcz}For every non-negative integers $n,$ the
Binet's formulas for Tribonacci and Tribonacci-Lucas numbers are given as%
\begin{eqnarray*}
T_{-n} &=&\frac{\alpha ^{-n+1}}{\left( \alpha -\gamma \right) \left( \alpha
-\beta \right) }+\frac{\beta ^{-n+1}}{\left( \beta -\gamma \right) \left(
\beta -\alpha \right) }+\frac{\gamma ^{-n+1}}{\left( \gamma -\beta \right)
\left( \gamma -\alpha \right) }, \\
K_{-n} &=&\alpha ^{-n}+\beta ^{-n}+\gamma ^{-n}.
\end{eqnarray*}
\end{corollary}

\textit{Proof.} From Theorem \ref{theorem:sdfvbaxdreqwsax}, we have 
\begin{eqnarray*}
\mathcal{T}_{-n} &=&\frac{\alpha ^{-n}}{\left( \alpha -\gamma \right) \left(
\alpha -\beta \right) }\left( 
\begin{array}{ccc}
\alpha ^{2} & \alpha ^{2}\left( \alpha -1\right) & \alpha \\ 
\alpha & \alpha ^{2}-\alpha & \alpha ^{2}\left( \alpha -1\right) -\alpha \\ 
\alpha ^{2}\left( \alpha -1\right) -\alpha & 2\alpha -\alpha ^{2}\left(
\alpha -1\right) & \alpha ^{2}-\alpha ^{2}\left( \alpha -1\right)%
\end{array}%
\right) \\
&&+\frac{\beta ^{-n}}{\left( \beta -\gamma \right) \left( \beta -\alpha
\right) }\left( 
\begin{array}{ccc}
\beta ^{2} & \beta ^{2}\left( \beta -1\right) & \beta \\ 
\beta & \beta ^{2}-\beta & \beta ^{2}\left( \beta -1\right) -\beta \\ 
\beta ^{2}\left( \beta -1\right) -\beta & 2\beta -\beta ^{2}\left( \beta
-1\right) & \beta ^{2}-\beta ^{2}\left( \beta -1\right)%
\end{array}%
\right) \\
&&+\frac{\gamma ^{-n}}{\left( \gamma -\beta \right) \left( \gamma -\alpha
\right) }\left( 
\begin{array}{ccc}
\gamma ^{2} & \gamma ^{2}\left( \gamma -1\right) & \gamma \\ 
\gamma & \gamma ^{2}-\gamma & \gamma ^{2}\left( \gamma -1\right) -\gamma \\ 
\gamma ^{2}\left( \gamma -1\right) -\gamma & 2\gamma -\gamma ^{2}\left(
\gamma -1\right) & \gamma ^{2}-\gamma ^{2}\left( \gamma -1\right)%
\end{array}%
\right)
\end{eqnarray*}%
By Theorem \ref{theorem:gbvdcsxfazoerstaqw}, we know that 
\begin{equation*}
\mathcal{T}_{-n}=\left( 
\begin{array}{ccc}
T_{-n+1} & T_{-n}+T_{-n-1} & T_{-n} \\ 
T_{-n} & T_{-n-1}+T_{-n-2} & T_{-n-1} \\ 
T_{-n-1} & T_{-n-2}+T_{-n-3} & T_{-n-2}%
\end{array}%
\right) .
\end{equation*}%
Now, if we compare the 2nd row and 1st column entries with the matrices in
the above two equations, then we obtain%
\begin{eqnarray*}
T_{-n} &=&\frac{\alpha ^{-n}\alpha }{\left( \alpha -\gamma \right) \left(
\alpha -\beta \right) }+\frac{\beta ^{-n}\beta }{\left( \beta -\gamma
\right) \left( \beta -\alpha \right) }+\frac{\gamma ^{-n}\gamma }{\left(
\gamma -\beta \right) \left( \gamma -\alpha \right) } \\
&=&\frac{\alpha ^{-n+1}}{\left( \alpha -\gamma \right) \left( \alpha -\beta
\right) }+\frac{\beta ^{-n+1}}{\left( \beta -\gamma \right) \left( \beta
-\alpha \right) }+\frac{\gamma ^{-n+1}}{\left( \gamma -\beta \right) \left(
\gamma -\alpha \right) }.
\end{eqnarray*}%
Tribonacci-Lucas case cen be proved similarly.

Now, we present summation formulas for Tribonacci and Tribonacci-Lucas
matrix sequences.

\begin{theorem}
For $m>j\geq 0,$ we have%
\begin{equation}
\sum_{i=0}^{n-1}\mathcal{T}_{-mi-j}=\frac{\mathcal{T}_{-mn+m-j}+\mathcal{T}%
_{-mn-m-j}+(1-K_{-m})\mathcal{T}_{-mn-j}}{K_{-m}-K_{m}}-\frac{\mathcal{T}%
_{-m-j}+\mathcal{T}_{-j+m}+(1-K_{-m})\mathcal{T}_{-j}}{K_{-m}-K_{m}}
\label{equati:tyujnbvcdfs}
\end{equation}%
and%
\begin{equation}
\sum_{i=0}^{n-1}\mathcal{K}_{-mi-j}=\frac{\mathcal{K}_{-mn+m-j}+\mathcal{K}%
_{-mn-m-j}+(1-K_{-m})\mathcal{K}_{-mn-j}}{K_{-m}-K_{m}}-\frac{\mathcal{K}%
_{-m-j}+\mathcal{K}_{-j+m}+(1-K_{-m})\mathcal{K}_{-j}}{K_{-m}-K_{m}}.
\label{equat:gbgdcvsxdzsarwt}
\end{equation}
\end{theorem}

\textit{Proof. Note that} 
\begin{eqnarray*}
\sum_{i=0}^{n-1}\mathcal{T}_{-mi-j} &=&\sum_{i=0}^{n-1}(A_{3}\alpha
^{-mi-j}+B_{3}\beta ^{-mi-j}+C_{3}\gamma ^{-mi-j}) \\
&=&A_{3}\alpha ^{-j}\left( \frac{\alpha ^{-mn}-1}{\alpha ^{-m}-1}\right)
+B_{3}\beta ^{-j}\left( \frac{\beta ^{-mn}-1}{\beta ^{-m}-1}\right)
+C_{3}\gamma ^{-j}\left( \frac{\gamma ^{-mn}-1}{\gamma ^{-m}-1}\right)
\end{eqnarray*}%
and%
\begin{eqnarray*}
\sum_{i=0}^{n-1}\mathcal{K}_{-mi-j} &=&\sum_{i=0}^{n-1}(A_{2}\alpha
^{-mi-j}+B_{2}\beta ^{-mi-j}+C_{2}\gamma ^{-mi-j}) \\
&=&A_{2}\alpha ^{-j}\left( \frac{\alpha ^{-mn}-1}{\alpha ^{-m}-1}\right)
+B_{2}\beta ^{-j}\left( \frac{\beta ^{-mn}-1}{\beta ^{-m}-1}\right)
+C_{2}\gamma ^{-j}\left( \frac{\gamma ^{-mn}-1}{\gamma ^{-m}-1}\right) .
\end{eqnarray*}%
Simplifying and rearranging the last equalities in the last two expression
imply (\ref{equati:tyujnbvcdfs}) and (\ref{equat:gbgdcvsxdzsarwt}) as
required. 
\endproof%

As in Corollary \ref{corollary:fvxgfsypqavbxcz}, in the proof of next
Corollary, we just compare the linear combination of the 2nd row and 1st
column entries of the relevant matrices.

\begin{corollary}
For $m>j>0,$ we have%
\begin{equation}
\sum_{i=0}^{n-1}T_{-mi-j}=\frac{T_{-mn+m-j}+T_{-mn-m-j}+(1-K_{-m})T_{-mn-j}}{%
K_{-m}-K_{m}}-\frac{T_{-m-j}+T_{-j+m}+(1-K_{-m})T_{-j}}{K_{-m}-K_{m}}
\end{equation}%
and%
\begin{equation}
\sum\limits_{i=0}^{n-1}K_{-mi-j}=\frac{%
K_{-mn+m-j}+K_{-mn-m-j}+(1-K_{-m})K_{-mn-j}}{K_{-m}-K_{m}}-\frac{%
K_{-m+j}+K_{-j+m}+(1-K_{-m})K_{-j}}{K_{-m}-K_{m}}
\end{equation}
\end{corollary}

Note that using the above Corollary we obtain the following well known
formulas (taking $m=1,j=0$):%
\begin{eqnarray*}
\sum_{i=0}^{n-1}T_{-i} &=&\frac{T_{-n+1}+T_{-n-1}+2T_{-n}-1}{-2}, \\
\sum_{i=0}^{n-1}K_{-i} &=&\frac{K_{-n+1}+K_{-n-1}+2K_{-n}-6}{-2}.\text{ \ }
\end{eqnarray*}

We now give generating functions of $\mathcal{T}$ and $\mathcal{K}$ for
negative indices.

\begin{theorem}
\label{theorem:amnbvcsertopu}For negative indices, the generating function
for the Tribonacci and Tribonacci-Lucas matrix sequences are given as%
\begin{equation*}
\sum_{n=0}^{\infty }\mathcal{T}_{-n}x^{n}=\frac{1}{1+x+x^{2}-x^{3}}\left( 
\begin{array}{ccc}
x^{2}+x+1 & x^{2}+x & x^{2} \\ 
x^{2} & x+1 & x \\ 
x & x^{2}-x & 1%
\end{array}%
\right)
\end{equation*}%
and%
\begin{equation*}
\sum_{n=0}^{\infty }\mathcal{K}_{-n}x^{n}=\frac{1}{1+x+x^{2}-x^{3}}\left( 
\begin{array}{ccc}
3x^{2}+4x+1 & 4x^{2}+2 & x^{2}+2x+3 \\ 
x^{2}+2x+3 & 2x^{2}+2x-2 & 3x^{2}-2x-1 \\ 
3x^{2}-2x-1 & -2x^{2}+4x+4 & -x^{2}+4x-1%
\end{array}%
\right)
\end{equation*}%
respectively.
\end{theorem}

\textit{Proof. }Then, using Definition \ref{definition::negabvhsyteur}, and
adding $xg(x)$ and $x^{2}g(x)$ to $g(x)$ and also substracting $x^{3}g(x)$
we obtain (note the shift in the index $n$ in the third line) 
\begin{eqnarray*}
(1+x+x^{2}-x^{3})g(x) &=&\sum_{n=0}^{\infty }\mathcal{T}_{-n}x^{n}+x%
\sum_{n=0}^{\infty }\mathcal{T}_{-n}x^{n}+x^{2}\sum_{n=0}^{\infty }\mathcal{T%
}_{-n}x^{n}-x^{3}\sum_{n=0}^{\infty }\mathcal{T}_{-n}x^{n} \\
&=&\sum_{n=0}^{\infty }\mathcal{T}_{-n}x^{n}+\sum_{n=0}^{\infty }\mathcal{T}%
_{-n}x^{n+1}+\sum_{n=0}^{\infty }\mathcal{T}_{-n}x^{n+2}-\sum_{n=0}^{\infty }%
\mathcal{T}_{-n}x^{n+3} \\
&=&\sum_{n=0}^{\infty }\mathcal{T}_{-n}x^{n}+\sum_{n=1}^{\infty }\mathcal{T}%
_{-n+1}x^{n}+\sum_{n=2}^{\infty }\mathcal{T}_{-n+2}x^{n}-\sum_{n=3}^{\infty }%
\mathcal{T}_{-n+3}x^{n} \\
&=&(\mathcal{T}_{0}+\mathcal{T}_{-1}x+\mathcal{T}_{-2}x^{2})+(\mathcal{T}%
_{0}x+\mathcal{T}_{-1}x^{2})+\mathcal{T}_{0}x^{2} \\
&&+\sum_{n=3}^{\infty }(\mathcal{T}_{-n}+\mathcal{T}_{-n+1}+\mathcal{T}%
_{-n+2}-\mathcal{T}_{-n+3})x^{n} \\
&=&(\mathcal{T}_{0}+\mathcal{T}_{-1}x+\mathcal{T}_{-2}x^{2})+(\mathcal{T}%
_{0}x+\mathcal{T}_{-1}x^{2})+\mathcal{T}_{0}x^{2} \\
&=&\mathcal{T}_{0}+(\mathcal{T}_{0}+\mathcal{T}_{-1})x+(\mathcal{T}_{0}+%
\mathcal{T}_{-1}+\mathcal{T}_{-2})x^{2}
\end{eqnarray*}%
Rearranging above equation, we get%
\begin{equation*}
g(x)=\frac{\mathcal{T}_{0}+(\mathcal{T}_{0}+\mathcal{T}_{-1})x+(\mathcal{T}%
_{0}+\mathcal{T}_{-1}+\mathcal{T}_{-2})x^{2}}{1+x+x^{2}-x^{3}}
\end{equation*}%
which equals the $\sum_{n=0}^{\infty }\mathcal{T}_{n}x^{n}$ in the Theorem.

Tribonacci-Lucas case can be proved similarly.%
\endproof%

Now, we will obtain generating functions for Tribonacci and Tribonacci-Lucas
numbers in terms of Tribonacci and Tribonacci-Lucas matrix sequences with
negative indices as a consequence of Theorem \ref{theorem:amnbvcsertopu}. To
do this, we will again compare the the 2nd row and 1st column entries with
the matrices in Theorem \ref{theorem:amnbvcsertopu}. Thus we have the
following corollary.

\begin{corollary}
\bigskip The generating functions for the Tribonacci sequence $%
\{T_{-n}\}_{n\geq 0}$ and Tribonacci-Lucas sequence $\{K_{-n}\}_{n\geq 0}$
are given as%
\begin{equation*}
\sum_{n=0}^{\infty }T_{-n}x^{n}=\frac{x^{2}}{1+x+x^{2}-x^{3}}\text{\ and \ }%
\sum_{n=0}^{\infty }K_{-n}x^{n}=\frac{x^{2}+2x+3}{1+x+x^{2}-x^{3}}.
\end{equation*}%
respectively.
\end{corollary}

Note that using above Corollary we can obtain Binet formulas for $T_{-n}$
and $K_{-n}$ again.

\section{Relation Between Tribonacci and Tribonacci-Lucas Matrix Sequences
With Negative Indices}

The following theorem shows that there always exist interrelation between
Tribonacci and Tribonacci-Lucas matrix sequences with negative indices.

\begin{theorem}
\label{theorem:dfgbvsopeqswadz}For the matrix sequences $\{\mathcal{T}%
_{-n}\}_{n\geq 0}$ and $\{\mathcal{K}_{-n}\}_{n\geq 0},\ $we have the
following identities.

\begin{description}
\item[(a)] $\mathcal{K}_{-n}=3\mathcal{T}_{-n+1}-2\mathcal{T}_{-n}-\mathcal{T%
}_{-n-1},$

\item[(b)] $\mathcal{K}_{-n}=\mathcal{T}_{-n}+2\mathcal{T}_{-n-1}+3\mathcal{T%
}_{-n-2},$

\item[(c)] $\mathcal{K}_{-n}=-\mathcal{T}_{-n+2}+4\mathcal{T}_{-n+1}-%
\mathcal{T}_{-n},$

\item[(d)] $\mathcal{T}_{-n}=\frac{1}{22}(5\mathcal{K}_{-n+2}-3\mathcal{K}%
_{-n+1}-4\mathcal{K}_{-n})$
\end{description}
\end{theorem}

\textit{Proof.} From (\ref{equation:gfdvbxczvsadou}), (\ref%
{equation:fvcdxsartewqa}) and (\ref{equation:yuhtgsdafcscvb}), (a), (b) (c)
follow. It is easy to show that $K_{-n}=-T_{-n+2}+4T_{-n+1}-T_{-n}$ and $%
22T_{-n}=5K_{-n+2}-3K_{-n+1}-4K_{-n},$ so now (d) and (e) follow. 
\endproof%

\begin{lemma}
\label{lemma:bcsertbvcxsaweqa}For all non-negative integers $m$ and $n,\ $we
have the following identities.

\begin{description}
\item[(a)] $\mathcal{K}_{0}\mathcal{T}_{-n}=\mathcal{T}_{-n}\mathcal{K}_{0}=%
\mathcal{K}_{-n},$

\item[(b)] $\mathcal{T}_{0}\mathcal{K}_{-n}=\mathcal{K}_{-n}\mathcal{T}_{0}=%
\mathcal{K}_{-n}.$
\end{description}
\end{lemma}

\textit{Proof.} Identities can be established easily. Note that to show (a)
we need to use all the relations (\ref{equation:gfdvbxczvsadou}), (\ref%
{equation:fvcdxsartewqa}) and (\ref{equation:yuhtgsdafcscvb}). 
\endproof%

Next Corollary gives another relation between the numbers $T_{-n}$ and $%
K_{-n}$ and also the matrices $\mathcal{T}_{-n}$ and $\mathcal{K}_{-n}$.

\begin{corollary}
We have the following identities.

\begin{description}
\item[(a)] $T_{-n}=\frac{1}{22}(K_{-n}+5K_{-n-1}+2K_{-n+1}),$

\item[(b)] $\mathcal{T}_{-n}=\frac{1}{22}(\mathcal{K}_{-n}+5\mathcal{K}%
_{-n-1}+2\mathcal{K}_{-n+1}).$
\end{description}
\end{corollary}

\textit{Proof.} From Lemma \ref{lemma:bcsertbvcxsaweqa} (a), we know that $%
\mathcal{K}_{0}\mathcal{T}_{-n}=\mathcal{K}_{-n}.$ To show (a), use Theorem %
\ref{theorem:gbvdcsxfazoerstaqw} for the matrix $\mathcal{T}_{-n}$ and
calculate the matrix operation $\mathcal{K}_{0}^{-1}\mathcal{K}_{-n}$ and
then compare the 2nd row and 1st column entries with the matrices $\mathcal{T%
}_{-n}$ and $\mathcal{K}_{0}^{-1}\mathcal{K}_{-n}.$ Now (b) follows from
(a). 
\endproof%

The following theorem shows that there exist relation between the positive
indices and negative indices for Tribonacci matrix sequences.

\begin{theorem}
\label{theorem:sdfgbvcxsdferat}For $n\geq 0,\ $we have the following
identity:%
\begin{equation}
\mathcal{T}_{-n}=(\mathcal{T}_{n})^{-1}.  \label{equat:sdfghbnxctyqua}
\end{equation}
\end{theorem}

\textit{Proof.}\ We prove by mathematical induction. If $n=0$ then we have 
\begin{equation*}
\mathcal{T}_{0}=\left( 
\begin{array}{ccc}
1 & 0 & 0 \\ 
0 & 1 & 0 \\ 
0 & 0 & 1%
\end{array}%
\right) =\left( 
\begin{array}{ccc}
1 & 0 & 0 \\ 
0 & 1 & 0 \\ 
0 & 0 & 1%
\end{array}%
\right) ^{-1}=(\mathcal{T}_{0})^{-1}
\end{equation*}%
which is true and 
\begin{equation*}
\mathcal{T}_{-1}=\left( 
\begin{array}{ccc}
0 & 1 & 0 \\ 
0 & 0 & 1 \\ 
1 & -1 & -1%
\end{array}%
\right) =\left( 
\begin{array}{ccc}
1 & 1 & 1 \\ 
1 & 0 & 0 \\ 
0 & 1 & 0%
\end{array}%
\right) ^{-1}=(\mathcal{T}_{1})^{-1}
\end{equation*}%
which is true. Assume that the equality holds for $n\leq k.$ For $n=k+1,$ by
using Teorem \ref{theorem:olscdxfczdarepwout}, we obtain%
\begin{eqnarray*}
(\mathcal{T}_{k+1})^{-1} &=&(\mathcal{T}_{k}\mathcal{T}_{1})^{-1}=(\mathcal{T%
}_{1})^{-1}(\mathcal{T}_{k})^{-1}=\mathcal{T}_{-1}\mathcal{T}_{-k} \\
&=&\left( 
\begin{array}{ccc}
0 & 1 & 0 \\ 
0 & 0 & 1 \\ 
1 & -1 & -1%
\end{array}%
\right) \left( 
\begin{array}{ccc}
T_{-k+1} & T_{-k}+T_{-k-1} & T_{-k} \\ 
T_{-k} & T_{-k-1}+T_{-k-2} & T_{-k-1} \\ 
T_{-k-1} & T_{-k-2}+T_{-k-3} & T_{-k-2}%
\end{array}%
\right) \\
&=&\left( 
\begin{array}{ccc}
T_{-k} & T_{-k-1}+T_{-k-2} & T_{-k-1} \\ 
T_{-k-1} & T_{-k-2}+T_{-k-3} & T_{-k-2} \\ 
T_{1-k}-T_{-k-1}-T_{-k} & T_{-k}-2T_{-k-2}-T_{-k-3} & 
T_{-k}-T_{-k-1}-T_{-k-2}%
\end{array}%
\right) \\
&=&\left( 
\begin{array}{ccc}
T_{-k} & T_{-k-1}+T_{-k-2} & T_{-k-1} \\ 
T_{-k-1} & T_{-k-2}+T_{-k-3} & T_{-k-2} \\ 
T_{-k-2} & T_{-k-3}+T_{-k-4} & T_{-k-3}%
\end{array}%
\right) \\
&=&\left( 
\begin{array}{ccc}
T_{-(k+1)+1} & T_{-(k+1)}+T_{-(k+1)-1} & T_{-(k+1)} \\ 
T_{-(k+1)} & T_{-(k+1)-1}+T_{-(k+1)-2} & T_{-(k+1)-1} \\ 
T_{-(k+1)-1} & T_{-(k+1)-2}+T_{-(k+1)-3} & T_{-(k+1)-2}%
\end{array}%
\right) \\
&=&\mathcal{T}_{-(k+1)}
\end{eqnarray*}%
Thus, by induction on $n,$ this proves (\ref{equat:sdfghbnxctyqua}). 
\endproof%

In the following Theorem we will use the next Lemma.

\begin{lemma}[{[\protect\ref{soykan2018triposmatrix}]}]
\label{lemma:gbnhytrdvopxdser}Let $A_{1},B_{1},C_{1};A_{2},B_{2},C_{2}$ as
in Theorem \ref{theorem:sacasdmnbhgvc}. Then the following relations hold:%
\begin{eqnarray*}
A_{1}^{2} &=&A_{1},\text{ }B_{1}^{2}=B_{1},\text{ }C_{1}^{2}=C_{1}, \\
A_{1}B_{1} &=&B_{1}A_{1}=A_{1}C_{1}=C_{1}A_{1}=C_{1}B_{1}=B_{1}C_{1}=\left(
0\right) .
\end{eqnarray*}
\end{lemma}

\textit{\ \ \ \ \ \ As the following theorem shows }$\mathcal{T}_{-m}$ and $%
\mathcal{K}_{-m}$ has nice properties.

\begin{theorem}
\label{theorem:mbosunmvdc}For $n,m\geq 0,\ $we have the following identities:

\begin{description}
\item[(a)] $\mathcal{T}_{-m}\mathcal{T}_{-n}=\mathcal{T}_{-m-n}$

\item[(b)] $\mathcal{T}_{m}\mathcal{T}_{-n}=\mathcal{T}_{m-n}$

\item[(c)] $\mathcal{T}_{-m}\mathcal{K}_{-n}=\mathcal{K}_{-n}\mathcal{T}%
_{-m}=\mathcal{K}_{-m-n}$

\item[(d)] $\mathcal{T}_{m}\mathcal{K}_{-n}=\mathcal{K}_{m}\mathcal{T}_{-n}=%
\mathcal{K}_{m-n}$

\item[(e)] $\mathcal{K}_{-m}\mathcal{K}_{-n}=\mathcal{K}_{0}\mathcal{K}%
_{-m-n}$

\item[(f)] $\mathcal{K}_{m}\mathcal{K}_{-n}=\mathcal{K}_{0}\mathcal{K}_{m-n}$

\item[(g)] $\sum_{i=0}^{n-1}\mathcal{T}_{-mi-j}=(\mathcal{T}_{-mn+m-j}-%
\mathcal{T}_{m-j})(\mathcal{T}_{0}-\mathcal{T}_{m})^{-1}.$
\end{description}
\end{theorem}

\textit{Proof.} \ 

\begin{description}
\item[(a)] Using Theorems \ref{theorem:sdfgbvcxsdferat} and \ref%
{theorem:olscdxfczdarepwout},\ we obtain%
\begin{equation*}
\mathcal{T}_{-m}\mathcal{T}_{-n}=(\mathcal{T}_{m})^{-1}(\mathcal{T}%
_{n})^{-1}=(\mathcal{T}_{n}\mathcal{T}_{m})^{-1}=(\mathcal{T}_{m+n})^{-1}=%
\mathcal{T}_{-m-n}.
\end{equation*}

\item[(b)] Using (\ref{equation:hbvcdfrgsxzasewq}), Lemma \ref%
{lemma:gbnhytrdvopxdser} and Theorem \ref{theorem:olscdxfczdarepwout} we
obtain 
\begin{eqnarray*}
\mathcal{T}_{m}\mathcal{T}_{-n} &=&(A_{1}\alpha ^{m}+B_{1}\beta
^{m}+C_{1}\gamma ^{-m})(A_{3}\alpha ^{-n}+B_{3}\beta ^{n}+C_{3}\gamma ^{-n})
\\
&=&A_{1}^{2}\alpha ^{m-n}+B_{1}^{2}\beta ^{m-n}+C_{1}^{2}\gamma
^{m-n}+A_{1}\allowbreak B_{1}\alpha ^{m}\beta ^{-n}+B_{1}A_{1}\alpha
^{-n}\beta ^{m} \\
&&+A_{1}C_{1}\alpha ^{m}\gamma ^{-n}+C_{1}A_{1}\alpha ^{-n}\gamma
^{m}+B_{1}C_{1}\beta ^{m}\allowbreak \gamma ^{-n}+C_{1}B_{1}\beta
^{-n}\gamma ^{m} \\
&=&A_{1}\alpha ^{m-n}+B_{1}\beta ^{m-n}+C_{1}\gamma ^{m-n} \\
&=&\mathcal{T}_{m-n}.
\end{eqnarray*}

\item[(c)] By Lemma \ref{lemma:bcsertbvcxsaweqa} (a) we have%
\begin{equation*}
\mathcal{T}_{-m}\mathcal{K}_{-n}=\mathcal{T}_{-m}\mathcal{T}_{-n}\mathcal{K}%
_{0}.
\end{equation*}%
By (a) and again by Lemma \ref{lemma:bcsertbvcxsaweqa} (a) we obtain%
\begin{equation*}
\mathcal{T}_{-m}\mathcal{K}_{-n}=\mathcal{T}_{-m-n}\mathcal{K}_{0}=\mathcal{K%
}_{-m-n}.
\end{equation*}%
The other equality can be obtained similarly.

\item[(d)] By Lemma \ref{lemma:bcsertbvcxsaweqa} (a) we have%
\begin{equation*}
\mathcal{T}_{m}\mathcal{K}_{-n}=\mathcal{T}_{m}\mathcal{T}_{-n}\mathcal{K}%
_{0}.
\end{equation*}%
By (b) and again by Lemma \ref{lemma:bcsertbvcxsaweqa} (a) we obtain%
\begin{equation*}
\mathcal{T}_{m}\mathcal{K}_{-n}=\mathcal{T}_{m-n}\mathcal{K}_{0}=\mathcal{K}%
_{m-n}.
\end{equation*}%
The other equality can be obtained similarly.

\item[(e)] By Lemma \ref{lemma:bcsertbvcxsaweqa} (a) we have%
\begin{equation*}
\mathcal{K}_{-m}\mathcal{K}_{-n}=\mathcal{K}_{0}\mathcal{T}_{-m}\mathcal{K}%
_{-n}.
\end{equation*}%
By (c) we obtain%
\begin{equation*}
\mathcal{K}_{-m}\mathcal{K}_{-n}=\mathcal{K}_{0}\mathcal{K}_{-m-n}.
\end{equation*}%
The other equality can be obtained similarly.

\item[(f)] By Theorem \ref{theorem:olscdxfczdarepwout} we can write%
\begin{equation*}
\mathcal{K}_{m}\mathcal{K}_{-n}=\mathcal{K}_{0}\mathcal{T}_{m}\mathcal{K}%
_{-n}
\end{equation*}%
From (d), we obtain%
\begin{equation*}
\mathcal{K}_{m}\mathcal{K}_{-n}=\mathcal{K}_{0}\mathcal{K}_{m-n}.
\end{equation*}

\item[(g)] We use (b). Since 
\begin{equation*}
\left( \sum_{i=0}^{n-1}\mathcal{T}_{-mi-j}\right) \mathcal{T}%
_{m}=\sum_{i=0}^{n-1}\mathcal{T}_{-mi-j+m}=\mathcal{T}_{m-j}+\left(
\sum_{i=0}^{n-1}\mathcal{T}_{-mi-j}\right) -\mathcal{T}_{-m(n-1)-j},
\end{equation*}%
we obtain%
\begin{equation*}
\mathcal{T}_{-m(n-1)-j}-\mathcal{T}_{m-j}=\left( \sum_{i=0}^{n-1}\mathcal{T}%
_{-mi-j}\right) -\left( \sum_{i=0}^{n-1}\mathcal{T}_{-mi-j}\right) \mathcal{T%
}_{m}=\left( \sum_{i=0}^{n-1}\mathcal{T}_{-mi-j}\right) \left( \mathcal{T}%
_{0}-\mathcal{T}_{m}\right)
\end{equation*}%
and so%
\begin{equation*}
\sum_{i=0}^{n-1}\mathcal{T}_{-mi-j}=(\mathcal{T}_{-mn+m-j}-\mathcal{T}%
_{m-j})(\mathcal{T}_{0}-\mathcal{T}_{m})^{-1}.
\end{equation*}
\endproof%
\end{description}

\begin{theorem}
For all non-negative integers $m$ and $n,\ $we have the following identities.

\begin{description}
\item[(a)] $\mathcal{K}_{-m}\mathcal{K}_{-n}=\mathcal{K}_{-n}\mathcal{K}%
_{-m}=9\mathcal{T}_{-m-n+2}-12\mathcal{T}_{-m-n+1}-2\mathcal{T}_{-m-n}+4%
\mathcal{T}_{-m-n-1}+\mathcal{T}_{-m-n-2},$

\item[(b)] $\mathcal{K}_{-m}\mathcal{K}_{-n}=\mathcal{K}_{-n}\mathcal{K}%
_{-m}=\mathcal{T}_{-m-n}+4\mathcal{T}_{-m-n-1}+10\mathcal{T}_{-m-n-2}+12%
\mathcal{T}_{-m-n-3}+\allowbreak 9\mathcal{T}_{-m-n-4},$

\item[(c)] $\mathcal{K}_{-m}\mathcal{K}_{-n}=\mathcal{K}_{-n}\mathcal{K}%
_{-m}=\mathcal{T}_{-m-n}-8\mathcal{T}_{-m-n+1}+18\mathcal{T}_{-m-n+2}-8%
\mathcal{T}_{-m-n+3}+\mathcal{T}_{-m-n+4}.$
\end{description}
\end{theorem}

\textit{Proof.} \ \ \ 

\begin{description}
\item[(a)] Using Theorem \ref{theorem:mbosunmvdc} (a) and Theorem \ref%
{theorem:dfgbvsopeqswadz} (a) we obtain%
\begin{eqnarray*}
\mathcal{K}_{-m}\mathcal{K}_{-n} &=&(3\mathcal{T}_{-m+1}-2\mathcal{T}_{-m}-%
\mathcal{T}_{-m-1})(3\mathcal{T}_{-n+1}-2\mathcal{T}_{-n}-\mathcal{T}_{-n-1})
\\
&=&2\mathcal{T}_{-n}\mathcal{T}_{-m-1}-6\mathcal{T}_{-n}\mathcal{T}_{-m+1}+2%
\mathcal{T}_{-m}\mathcal{T}_{-n-1}-6\mathcal{T}_{-m}\mathcal{T}_{-n+1} \\
&&+\allowbreak 4\mathcal{T}_{-m}\mathcal{T}_{-n}+\mathcal{T}_{-m-1}\mathcal{T%
}_{-n-1}-3\mathcal{T}_{-m-1}\mathcal{T}_{-n+1}-3\mathcal{T}_{-m+1}\mathcal{T}%
_{-n-1}+\allowbreak 9\mathcal{T}_{-m+1}\mathcal{T}_{-n+1} \\
&=&2\mathcal{T}_{-m-n-1}-6\mathcal{T}_{-m-n+1}+2\mathcal{T}_{-m-n-1}-6%
\mathcal{T}_{-m-n+1}+\allowbreak 4\mathcal{T}_{-m-n}+\mathcal{T}_{-m-n-2}-3%
\mathcal{T}_{-m-n} \\
&&-3\mathcal{T}_{-m+n}+\allowbreak 9\mathcal{T}_{-m-n+2} \\
&=&9\mathcal{T}_{-m-n+2}-12\mathcal{T}_{-m-n+1}-2\mathcal{T}_{-m-n}+4%
\mathcal{T}_{-m-n-1}+\mathcal{T}_{-m-n-2}
\end{eqnarray*}
$\allowbreak $

It can be shown similarly that $\mathcal{K}_{-n}\mathcal{K}_{-m}=9\mathcal{T}%
_{-m-n+2}-12\mathcal{T}_{-m-n+1}-2\mathcal{T}_{-m-n}+4\mathcal{T}_{-m-n-1}+%
\mathcal{T}_{-m-n-2}.$

The remaining of identities can be proved by considering again (a) and
Theorem \ref{theorem:dfgbvsopeqswadz}. 
\endproof%
\end{description}

Comparing matrix entries and using Teorem \ref{theorem:gbvdcsxfazoerstaqw}
we have next two result.

\begin{corollary}
For Tribonacci and Tribonacci-Lucas numbers, we have the following
identities:

\begin{description}
\item[(a)] $T_{-m-n}=T_{-m}T_{-n+1}+T_{-n}\left( T_{-m-1}+T_{-m-2}\right)
+T_{-m-1}T_{-n-1}$

\item[(b)] $K_{-m-n}=T_{-m}K_{-n+1}+K_{-n}\left( T_{-m-1}+T_{-m-2}\right)
+K_{-n-1}T_{-m-1}$

\item[(c)] $K_{-m}K_{-n+1}+K_{-n}\left( K_{-m-1}+K_{-m-2}\right)
+K_{-m-1}K_{-n-1}=9T_{-m-n+2}-12T_{-m-n+1}-2T_{-m-n}+4T_{-m-n-1}+T_{-m-n-2}$

\item[(d)] $K_{-m}K_{-n+1}+K_{-n}\left( K_{-m-1}+K_{-m-2}\right)
+K_{-m-1}K_{-n-1}=\allowbreak
T_{-m-n}+4T_{-m-n-1}+10T_{-m-n-2}+12T_{-m-n-3}+9T_{-m-n-4}$

\item[(e)] $K_{-m}K_{-n+1}+K_{-n}\left( K_{-m-1}+K_{-m-2}\right)
+K_{-m-1}K_{-n-1}=T_{-m-n}-8T_{-m-n+1}+18T_{-m-n+2}-8T_{-m-n+3}+T_{-m-n+4}$
\end{description}
\end{corollary}

The following theorem shows that there exist relation between the positive
indices and negative indices for Tribonacci-Lucas matrix sequences.

\begin{theorem}
\label{theorem:hgfdscvzcxsad}For all non-negative integers $n,\ $we have the
following identity:%
\begin{equation*}
\mathcal{K}_{-n}=(\mathcal{K}_{0})^{1-n}(\mathcal{K}_{-1})^{n}
\end{equation*}
\end{theorem}

\textit{Proof.} \ Taking $(n-1)$ for $m$ and $1$ for $n$ in $\mathcal{K}_{-m}%
\mathcal{K}_{-n}=\mathcal{K}_{0}\mathcal{K}_{-m-n}$ which is given in
Theorem \ref{theorem:mbosunmvdc} (e), we obtain that%
\begin{equation}
\mathcal{K}_{0}\mathcal{K}_{-n}=\mathcal{K}_{-n+1}\mathcal{K}_{-1}.
\label{equati:nvcxsdxzcaewq}
\end{equation}%
If we multiply both side of the equation (\ref{equati:nvcxsdxzcaewq}) with $%
\mathcal{K}_{0}$ we have the relation%
\begin{eqnarray*}
\mathcal{K}_{0}\mathcal{K}_{0}\mathcal{K}_{-n} &=&\mathcal{K}_{0}\mathcal{K}%
_{-n+1}\mathcal{K}_{-1} \\
&=&\mathcal{K}_{-n+2}\mathcal{K}_{-1}\mathcal{K}_{-1}.
\end{eqnarray*}%
Repeating this process we then obtain%
\begin{equation*}
\mathcal{K}_{0}^{n-1}\mathcal{K}_{-n}=\mathcal{K}_{-1}^{n}.
\end{equation*}%
Thus, it follows that%
\begin{equation*}
\mathcal{K}_{-n}=\mathcal{K}_{0}^{1-n}\mathcal{K}_{-1}^{n}.
\end{equation*}%
This completes the proof. 
\endproof%

Note that using Theorem \ref{theorem:mbosunmvdc} (d) in Theorem \ref%
{theorem:hgfdscvzcxsad}, we obtain%
\begin{equation*}
\mathcal{K}_{-n}=(\mathcal{K}_{n}\mathcal{T}_{-n})^{1-n}\mathcal{K}_{-1}^{n}=%
\mathcal{T}_{-n}^{1-n}\mathcal{K}_{n}^{1-n}\mathcal{K}_{-1}^{n}
\end{equation*}%
and then by Theorem \ref{theorem:sdfgbvcxsdferat} we get%
\begin{equation*}
\mathcal{K}_{-n}=(\mathcal{T}_{n}\mathcal{K}_{n}^{-1}\mathcal{K}_{-1})^{n-1}%
\mathcal{K}_{-1}.
\end{equation*}

The next two theorems provide us the convenience to obtain the powers of
Tribonacci and Tribonacci-Lucas matrix sequences.

\begin{theorem}
For non-negatif integers $m,n$ and $r$ with $n\geq r,$ the following
identities hold:

\begin{description}
\item[(a)] $(\mathcal{T}_{-n})^{m}=\mathcal{T}_{-mn},$

\item[(b)] $(\mathcal{T}_{-n-1})^{m}=(\mathcal{T}_{-1})^{m}\mathcal{T}%
_{-mn}, $

\item[(c)] $\mathcal{T}_{-n-r}\mathcal{T}_{-n+r}=(\mathcal{T}_{-n})^{2}=(%
\mathcal{T}_{-2})^{n}.$
\end{description}
\end{theorem}

\textit{Proof.} \ \ 

\begin{description}
\item[(a)] By Theorem \ref{theorem:sdfgbvcxsdferat} we have 
\begin{equation*}
(\mathcal{T}_{-n})^{m}=\mathcal{((T}_{n}\mathcal{)}^{-1}\mathcal{)}^{m}=%
\mathcal{((T}_{n}\mathcal{)}^{m}\mathcal{)}^{-1}.
\end{equation*}%
Using Theorem \ref{theorem:olscdxfczdarepwout}, and also\ again Theorem \ref%
{theorem:sdfgbvcxsdferat} we obtain%
\begin{equation*}
(\mathcal{T}_{-n})^{m}=(\mathcal{T}_{mn})^{-1}=\mathcal{T}_{-mn}.
\end{equation*}

\item[(b)] Using Theorem \ref{theorem:mbosunmvdc} (a) and method used in
above (a) we can write%
\begin{equation*}
(\mathcal{T}_{-n-1})^{m}=\mathcal{T}_{m(-n-1)}=\mathcal{T}_{-m}\mathcal{T}%
_{-mn}=\mathcal{T}_{-1}\mathcal{T}_{-m+1}\mathcal{T}_{-mn}.
\end{equation*}%
Similarly, we obtain $\mathcal{T}_{-m+1}=\mathcal{T}_{-1}\mathcal{T}_{-m+2}.$
Continuing to this iterative process, we then obtain%
\begin{equation*}
(\mathcal{T}_{-n-1})^{m}=\underset{m\text{ times}}{\underbrace{\mathcal{T}%
_{-1}\mathcal{T}_{-1}...\mathcal{T}_{-1}}}\mathcal{T}_{-mn}=(\mathcal{T}%
_{-1})^{m}\mathcal{T}_{-mn}.
\end{equation*}

\item[(c)] The proof is similar to (b). 
\endproof%
\end{description}

We have analogues results for the matrix sequence $\mathcal{K}_{-n}.$

\begin{theorem}
For non-negatif integers $m,n$ and $r$ with $n\geq r,$ the following
identities hold:

\begin{description}
\item[(a)] $\mathcal{K}_{-n-r}\mathcal{K}_{-n+r}=(\mathcal{K}_{-n})^{2},$

\item[(b)] $(\mathcal{K}_{-n})^{m}=\mathcal{K}_{0}^{m}\mathcal{T}_{-mn}=%
\mathcal{K}_{0}^{m-1}\mathcal{K}_{-mn}.$
\end{description}
\end{theorem}

\textit{Proof.} \ \ 

\begin{description}
\item[(a)] Applying Theorem \ref{theorem:mbosunmvdc} (e), we find%
\begin{equation*}
\mathcal{K}_{-n-r}\mathcal{K}_{-n+r}=\mathcal{K}_{0}\mathcal{K}_{-2n}=%
\mathcal{K}_{0}\mathcal{K}_{-n-n}=\mathcal{K}_{-n}\mathcal{K}_{-n}=(\mathcal{%
K}_{-n})^{2}.
\end{equation*}

\item[(b)] By Lemma \ref{lemma:bcsertbvcxsaweqa} (a), we see that%
\begin{equation*}
(\mathcal{K}_{-n})^{m}=(\mathcal{K}_{0}\mathcal{T}_{-n})^{m}=(\mathcal{K}%
_{0})^{m}(\mathcal{T}_{-n})^{m}.
\end{equation*}%
Using Theorem \ref{theorem:mbosunmvdc} (a), we obtain%
\begin{equation*}
(\mathcal{K}_{-n})^{m}=\mathcal{K}_{0}^{m}\mathcal{T}_{-mn}.
\end{equation*}%
Thus, we can write $(\mathcal{K}_{-n})^{m}=\mathcal{K}_{0}^{m-1}\mathcal{K}%
_{0}\mathcal{T}_{-mn}$ and by Lemma \ref{lemma:bcsertbvcxsaweqa} (a) we get $%
(\mathcal{K}_{-n})^{m}=\mathcal{K}_{0}^{m-1}\mathcal{K}_{-mn}.$ This
completes the proof. 
\endproof%
\end{description}

\section{On the Theorem of Rabinowitz and Bruckman}

We now present a remarkable theorem of Rabinowitz and Bruckman [\ref%
{bib:bruckman1995}].

\begin{theorem}
\label{theorem:dfgcvzoustraewqa}(Rabinowitz, 1994 [\ref{bib:rabinowitz1994}%
], solution: Bruckman, 1995 [\ref{bib:bruckman1995}]) Assume that $H_{n}$
satisfies a second-order linear recurrence with constant coefficients. Let $%
\{a_{i}\}$ and $\{b_{i}\},$ $i=1,2,...,r$ be integer constants and let $%
f(x_{0},x_{1},x_{2},...,x_{r})$ be a polynomial with integer coefficients.
If the expression 
\begin{equation*}
f((-1)^{n},H_{a_{1}n+b_{1}},H_{a_{2}n+b_{2}},...,H_{a_{r}n+b_{r}})
\end{equation*}%
vanishes for all integers $n>N$ then the expression vanishes for all
integers $n.$

(As a special case, if an identity involving Fibonacci and Lucas numbers is
true for all positive subscripts, then it must be true for all negative
subscripts as well.)
\end{theorem}

\textit{Proof. }
For details, see Bruckman [\ref{bib:bruckman1995}]. 
\endproof%

It follows from above Theorem that if an identity involving generalized
Fibonacci numbers (Horadam numbers) is true for all positive subscripts,
then it is true for all non-positive subscripts as well.

After seeing the results of this paper, we can propose the following
conjecture.

\begin{conjecture}
\label{theorem:mnbvyuhgasdfxcz}Assume that $H_{n}$ satisfies an order-k
linear homogeneous recurrences with constant coefficients. Let $\{a_{i}\}$
and $\{b_{i}\},$ $i=1,2,...,r$ be integer constants and let $%
f(x_{0},x_{1},x_{2},...,x_{r})$ be a polynomial with integer coefficients.
If the expression 
\begin{equation*}
f(H_{a_{1}n+b_{1}},H_{a_{2}n+b_{2}},...,H_{a_{r}n+b_{r}})
\end{equation*}%
vanishes for all integers $n>N$ then the expression vanishes for all
integers $n.$
\end{conjecture}

It seems that the proof of the Theorem \ref{theorem:dfgcvzoustraewqa} can be
applied to prove the above conjecture. It follows that if Conjecture \ref%
{theorem:mnbvyuhgasdfxcz} is true then if an identity involving generalized
Tribonacci numbers as a third-order linear recurrence with constant
coefficients (such as Tribonacci and Tribonacci-Lucas numbers) is true for
all positive subscripts, then it is true for all non-positive subscripts as
well.


\begin{thebibliography}{99}
\bibitem{akbulak2009} \label{akbulak2009}Akbulak, M., and Bozkurt, D., On
the Order-m Generalized Fibonacci k-numbers, Chaos Solitons\&Fractals,
42(3), p.1347-1355, 2009.

\bibitem{basu2014} \label{basu2014}Basu, M., and Das, M., Tribonacci
Matrices and a New Coding Theory, Discrete Mathematics, Algorithms and
Applications Vol. 06, No. 01 pp.1450008 (17 pages), 2014.

\bibitem{bruce1984} \label{bib:bruce1984}Bruce, I., A modified Tribonacci
sequence, The Fibonacci Quarterly, 22 : 3, pp. 244--246, 1984.

\bibitem{bruckman1995} \label{bib:bruckman1995}Paul S. Bruckman. P. S.,
"Solution to Problem H-487 (Proposed by Stanley Rabinowitz)." The Fibonacci
Quarterly 33, p.382, 1995.

\bibitem{bueno:2015} \label{bueno:2015}Bueno, A.C.F., A Note on Generalized
Tribonacci Sequence, Notes on Number Theory and Discrete Mathematics, Vol.
21, No. 1, 67--69, 2015.

\bibitem{cerdamoralez2018aa} \label{cerdamoralez2018aa}Cerda-Morales, G., On
the Third-Order Jabosthal and Third-Order Jabosthal-Lucas Sequences and
Their Matrix Representations, arxiv:1806.03709v1 [math.CO], 2018.

\bibitem{civciv2008} \label{civciv2008}Civciv, H., and Turkmen, R., On the
(s; t)-Fibonacci and Fibonacci matrix sequences, Ars Combin. 87, 161-173,
2008.

\bibitem{civciv2008b} \label{civciv2008b}Civciv, H., and Turkmen, R., Notes
on the (s; t)-Lucas and Lucas matrix sequences, Ars Combin. 89, 271-285,
2008.

\bibitem{cohen2002} \label{cohen2002}Cohen, M.S., Fernandez, J.J., Park,
A.E., and Schmedders, K., The Fibonacci Sequence: Relationship to the Human
Hand, The Journal of Hand Surgery, 28(1), p.157-160, 2002.

\bibitem{duchene2008} \label{duchene2008}Duch\^{e}ne, E., and Rigo, M., A
Morphic Approach to Combinatorial Games : the Tribonacci case, RAIRO -
Theoretical Informatics and Applications, Volume 42 (2008) no. 2 , p.
375-393.

\bibitem{feinberg1963} \label{bib:feinberg1963}Feinberg, M.,
Fibonacci--Tribonacci, The Fibonacci Quarterly, 1 : 3 (1963) pp. 71--74,
1963.

\bibitem{fiorenza2011} \label{fiorenza2011}Fiorenza, A., G. Vincenzi, G.,
Limit of ratio of consecutive terms for general order-k linear homogeneous
recurrences with constant coefficients, Chaos Solitons \& Fractals, 44(1-3)
p.147-152, 2011.

\bibitem{gulec2012} \label{gulec2012}Gulec, H.H., and Taskara, N., On the
(s; t)-Pell and (s; t)-Pell-Lucas sequences and their matrix
representations, Appl. Math. Lett. 25, 1554-1559, 2012.

\bibitem{howard2010} \label{bib:howard2010}Howard, F.T., Saidak, F.,
Congress Numer. 200 (2010), 225-237, 2010.

\bibitem{marohnic2012} \label{bib:marohnic2012}Marohni\'{c}, L., Strme\v{c}%
ki, T., Plastic Number: Construction and Applications, Advanced Research in
Scientific Areas 2012, 1523-1528, 2012.

\bibitem{rabinowitz1994} \label{bib:rabinowitz1994}Rabinowitz, S., Problem
H-487, The Fibonacci Quarterly 32, p.187, 1994.

\bibitem{randic1996} \label{randic1996}Randi\'{c}., M., Morales, D.A.,
Araujo, O., Higher-order Fibonacci Numbers, Journal of Mathematical
Chemistry, 20(1), p.79-94, 1996.

\bibitem{ridley1982} \label{ridley1982}Ridley, J.N., Packing Efficiency in
Sunflower Heads, Mathematical Biosciences, 58(1), p.129-139, 1982.

\bibitem{piezas} \label{piezas}Piezas, T., A Tale of Four Constants,
https://sites.google.com/site/tpiezas/0012.

\bibitem{scott1977} \label{bib:scott1977}Scott, A., Delaney, T., Hoggatt
Jr.,\ V., The Tribonacci sequence, The Fibonacci Quarterly, 15:3, pp.
193--200, 1977.

\bibitem{shannon1977} \label{bib:shannon1977}Shannon, A., Tribonacci numbers
and Pascal's pyramid, The Fibonacci Quarterly, 15:3, pp. 268-275, 1977.

\bibitem{sloane} \label{bib:sloane}N.J.A. Sloane, The on-line encyclopedia
of integer sequences, http://oeis.org/

\bibitem{soykan2018triposmatrix} \label{soykan2018triposmatrix}Soykan, Y.,
Matrix Sequences of Tribonacci and Tribonacci-Lucas Numbers,
arXiv:1809.07809v1 [math.NT]\ 20 Sep 2018.

\bibitem{spickerman1981} \label{bib:spickerman1981}Spickerman, W., Binet's
formula for the Tribonacci sequence, The Fibonacci Quarterly, 20,
pp.118--120, 1981.

\bibitem{uslu2013} \label{uslu2013}Uslu, K., and Uygun, S., On the (s,t)
Jacobsthal and (s,t) Jacobsthal-Lucas Matrix Sequences, Ars Combin. 108,
13-22, 2013.

\bibitem{uygun2015} \label{uygun2015}Uygun, \c{S}., and Uslu, K.,
(s,t)-Generalized Jacobsthal Matrix Sequences, Springer Proceedings in
Mathematics\&Statistics, Computational Analysis, Amat, Ankara, May 2015,
325-336.

\bibitem{uygun2016} \label{uygun2016}Uygun, \c{S}., Some Sum Formulas of
(s,t)-Jacobsthal and (s,t)-Jacobsthal Lucas Matrix Sequences, Applied
Mathematics, 7, 61-69, 2016.

\bibitem{yalavigi1972} \label{bib:yalavigi1972}Yalavigi, C. C., Properties
of Tribonacci numbers, The Fibonacci Quarterly, 10 : 3, pp. 231--246, 1972.

\bibitem{yazlik2013a} \label{yazlik2013a}Yazlik, Y., and Taskara, N., Uslu,
K. and Yilmaz, N. The generalized (s; t)-sequence and its matrix sequence,
Am. Inst. Phys. (AIP) Conf. Proc. 1389, 381-384, 2012.

\bibitem{yegnanarayanan2012} \label{yegnanarayanan2012}Yegnanarayanan V.,
The Chromatic number of Generalized Fibonacci Prime Distance Graph, Journal
of Math and Computatioal Science, 2 N.5 1451-1463, 2012.

\bibitem{yilmaz2013} \label{yilmaz2013}Yilmaz, N., and Taskara, N., Matrix
Sequences in Terms of Padovan and Perrin Numbers, Journal of Applied
Mathematics, Volume 2013, Article ID 941673, 7 pages, 2013.

\bibitem{yilmaz2014a} \label{yilmaz2014a}Yilmaz, N., Taskara, N., On the
Negatively Subscripted Padovan and Perrin Matrix Sequences, Communications
in Mathematics and Applications, Vol. 5, No. 2, 59-72, 2014.

\bibitem{wani2018} \label{wani2018}Wani, A.A., Badshah, V.H., and Rathore,
G.B.S., Generalized Fibonacci and k-Pell Matrix Sequences, Punjab University
Journal of Mathematics (ISSN 1016-2526) Vol. 50(1) (2018) pp. 68-79.
\end{thebibliography}
\end{document}